\documentclass[a4paper, reqno, 11pt]{amsart}
\usepackage[scale=0.72,centering]{geometry}
\usepackage{amsmath,mathabx,amssymb,amsthm,graphicx,bbm,pdfpages}
\usepackage{color,mathrsfs,tikz,hyperref}

\usepackage[font={footnotesize}]{caption}
\usepackage{amsaddr}

\numberwithin{equation}{section}

\newcommand{\bbN}{{\ensuremath{\mathbbm N}} }

\newcommand{\bbR}{{\ensuremath{\mathbbm R}} }

\newcommand{\bbZ}{{\ensuremath{\mathbbm Z}} }

\newcommand{\cC}{{\ensuremath{\mathcal C}} }

\newcommand{\cG}{{\ensuremath{\mathcal G}} }
\newcommand{\cH}{{\ensuremath{\mathcal H}} }

\newcommand{\bG}{{\ensuremath{\mathbf G}} }

\newtheorem{assump}{Assumption}
\newtheorem{problem}{Problem}
\newtheorem{theorem}{Theorem}[section]

\newtheorem{lemma}[theorem]{Lemma}
\newtheorem{corollary}[theorem]{Corollary}

\newtheorem{remark}{Remark}
\newtheorem{conjecture}{Conjecture}
\newtheorem{condition}{Condition}

\author[C.-H. Huang]{Chien-Hao Huang*}
\address[Chien-Hao Huang]{*Department of Mathematical Sciences\\
	National Chengchi University\\
	Taipei, Taiwan}
\address{*Corresponding author. E-mail: p092221010@gmail.com}

\email[C.~Huang]{p092221010@gmail.com}

\begin{document}

\title[Gaussian comparison for variance]{Nonsymmetric examples for Gaussian correlation inequalities}

\begin{abstract}
In this paper, we compare two variances of maxima of $N$ standard Gaussian random variables. One is a sequence of $N$ i.i.d. standard Gaussians, and the other one is $N$ standard Gaussians with covariances $\sigma_{1,2}=\rho \in(0,1)$ and $\sigma_{i,j}=0$, for other $ i\neq j$. It turns out that we need to discuss the covariance of two functions with respect to multivariate Gaussian distributions. Gaussian correlation inequalities hold for many symmetric (with respect to the origin) cases. However, in our case, the max function and its derivatives are not symmetric about the origin. We have two main results in this paper. First, we prove a specific case for a convex/log-concave correlation inequality for the standard multivariate Gaussian distribution. The other result is that the variance of maxima of standard Gaussians with $\sigma_{1,2}=\rho \in(0,1)$, while $\sigma_{i,j}=0$, for other $ i\neq j$, is larger than the variance of maxima of independent standard Gaussians. This implies that the variance of maxima of $N$ i.i.d. standard Gaussians is decreasing in $N$.
\end{abstract}


\maketitle


\section{Introduction}

We are interested in the following problem. Let $X$ and $Y$ be two random vectors with $X\sim Normal(0_N, \Sigma^X)$, $Y\sim Normal(0_N, \Sigma^Y)$.

\begin{assump}\label{assump1}
$N\geq 3$ and	$	\sigma^Y_{i,i} =\sigma^X_{i,i}=1, \;\; i=1,...,N$.
\end{assump}

Denote the max function as $M_N(x)=\max\limits_{i=1,...,N} x_i$. 

\begin{problem}\label{problem1}
	With Assumption \ref{assump1}, under what condition 
	\begin{equation}\text{Var}(M_N(Y)) \geq \text{Var}(M_N(X)) \;\; ?
	\end{equation}
\end{problem}

For the comparison between the expected values, Slepian \cite{Slepian62} proved the following.
\begin{theorem}(\cite{Slepian62}) Suppose that $\sigma^Y_{i,i} =\sigma^X_{i,i}, \;\; i=1,...,N$ and $E[Y_iY_j]\leq E[X_iX_j]$ for all $i,j$,
then
\begin{equation} P(M_N(Y) > u) \geq  P(M_N(X) > u) \end{equation}
for all real $u$. Moreover,	
	\begin{equation} E\left[ M_N(Y)\right]\geq E\left[ M_N(X)\right]. \end{equation}
\end{theorem}
See more details in Sect. 2.2, \cite{AT07}.\\

To attack Problem \ref{problem1}, we apply the classical method to the variance. We consider a smooth function $\phi:\bbR^N\to \bbR$ and
\begin{equation}\label{VarDiff_01}
	 \text{Var}(\phi(Y))- \text{Var}(\phi(X))=E\phi^2(Y)- E\phi^2(X) -\left[ (E\phi(Y))^2-(E\phi(X))^2 \right].
\end{equation}
\noindent
Later, we will use a specific $\phi$ to approach $M_N$. We now interpolate $X$ and $Y$. For any $0\leq \theta \leq 1$, let $Z(\theta)\sim Normal (0_N, \Sigma (\theta) )$,
where $\Sigma (\theta)=(1-\theta)\Sigma_X+\theta\Sigma_Y$, so that $Z(1)=Y$ and $Z(0)=X$.
Denote the probability density function of $Z(\theta)$ as $p_\theta(z)$.
Let $\hat{\psi}(\theta)=E\phi^2(Z(\theta))$ and $\psi(\theta)=(E\phi(Z(\theta)))^2$. \eqref{VarDiff_01} is equal to 
\begin{equation}\label{VarDiff_02}
	\hat{\psi}(1)-\hat{\psi}(0) -\left[\psi(1)-\psi(0)\right] =\int_0^1 \left\{\hat{\psi}'(\theta)- \psi'(\theta)\right\}d\theta .
\end{equation}

The first term in the integrand of \eqref{VarDiff_02} is 
\begin{equation}\label{hatphi1}
	\hat{\psi}'(\theta)= \frac{d}{d\theta} E\phi^2(Z(\theta)) = \int_{\bbR^N}  \phi^2(z) \:\frac{d}{d\theta} p_\theta(z) dz.
\end{equation}
It is known that the RHS of \eqref{hatphi1} is equal to
\begin{equation}\label{hatphi2}
	\int_{\bbR^N}  \phi^2(z) \:\frac{d}{d\theta} p_\theta(z) dz = \frac{1}{2} \sum_{i,j} \left(\sigma^Y_{i,j}-\sigma^X_{i,j}\right)\int_{\bbR^N}  \phi^2(z) \:\frac{\partial^2}{\partial z_i \partial z_j} p_\theta(z) dz.
\end{equation}
With integration by parts, \eqref{hatphi1} and \eqref{hatphi2} give
\begin{equation}\label{hatphi3}
	\hat{\psi}'(\theta)= \frac{1}{2} \sum_{i,j} \left(\sigma^Y_{i,j}-\sigma^X_{i,j}\right)\int_{\bbR^N}  \left[2\frac{\partial}{\partial x_j } \phi(z) \frac{\partial}{\partial x_i } \phi(z) +2\phi(z)\frac{\partial^2}{\partial x_j \partial x_i} \phi(z) \right]\:p_\theta(z) dz.
\end{equation}
Using the same way for $\psi'(\theta)$, we have
\begin{equation}\label{phi3} \begin{array}{rcl}
		\psi'(\theta)&=& \displaystyle 2 E(\phi(Z(\theta))) \cdot \frac{d}{d\theta} E(\phi(Z(\theta)))\\
		&=& \displaystyle 2 E(\phi(Z(\theta))) \cdot \frac{1}{2} \sum_{i,j} \left(\sigma^Y_{i,j}-\sigma^X_{i,j}\right)\int_{\bbR^N}  \left[ \frac{\partial^2}{\partial x_j \partial x_i} \phi(z) \right]\:p_\theta(z) dz \\
		&=& \displaystyle \frac{1}{2} \sum_{i,j}\limits \left(\sigma^Y_{i,j}-\sigma^X_{i,j}\right) 2 E(\phi(Z(\theta))) \cdot E\left[\frac{\partial^2}{\partial x_j \partial x_i}\phi(Z(\theta))\right].
	\end{array}
\end{equation}
\noindent
Summarizing \eqref{VarDiff_02}, \eqref{hatphi3} and \eqref{phi3}, \eqref{VarDiff_01} becomes
\begin{equation}\label{VarDiff_03}  \begin{split}
		&   \displaystyle \text{Var}(\phi(Y))- \text{Var}(\phi(X))\\
		=&  \displaystyle \frac{1}{2} \sum_{i,j} \left(\sigma^Y_{i,j}-\sigma^X_{i,j}\right)\\
		&\cdot \int_0^1 \;  
		\left\{ E\left[2\frac{\partial}{\partial x_j } \phi(Z(\theta)) \frac{\partial}{\partial x_i } \phi(Z(\theta))\right]
		+\text{Cov}\left(2\phi(Z(\theta)) , \frac{\partial^2}{\partial x_j \partial x_i} \phi(Z(\theta))\right) \right\}d\theta.
	\end{split}
\end{equation}

If one want to show that $\text{Var}(\phi(Y))- \text{Var}(\phi(X))\geq 0$ with the equality \eqref{VarDiff_03}
and the following condition,
\begin{condition}\label{condition1}
	$N\geq 3$ and $\sigma^Y_{i,j}\geq \sigma^X_{i,j},\; i\neq j$ with $\sigma^Y_{i,i}= \sigma^X_{i,i},\; \forall i,j =1,...,N$,
\end{condition}
\noindent
it then suffices to show that for each $i\neq j, \theta$ in the big parentheses in \eqref{VarDiff_03}, 
\begin{equation}\label{VarDiff_04}
	\frac{\partial \phi}{\partial x_j }  \frac{\partial \phi}{\partial x_i } \geq 0
	\;\; \text{and} \;\;\text{Cov}\left(\phi(Z(\theta)) , \frac{\partial^2}{\partial x_j \partial x_i} \phi(Z(\theta))\right)\geq 0.
\end{equation}

We go back to discuss Problem \ref{problem1}. For a constant $\beta>0$, we take $\phi(x)=Q_\beta (x):= \beta^{-1} \log S_N(x)$, where $S_N(x):= \sum\limits_{1\leq i \leq N} e^{\beta x_i}$. We plan to use $Q_\beta(x)$ to approximate $M_N (x)=\max\limits_{i=1,...,N} x_i$. Notice that

\begin{equation} \label{inequality_max_sum}
	M_N \leq Q_\beta \leq \frac{1}{\beta}\log N + M_N
\end{equation}
connects $Q_\beta$ and $M_N$. Moreover, $\lim\limits_{\beta\to \infty} Q_\beta =M_N$.

With simple calculations,
\begin{equation} \label{partial_max}
	\frac{\partial Q_\beta}{\partial x_{i}}=p_i:=\frac{e^{\beta x_i}}{S_N(x)},
	\;\; \frac{\partial p_i}{\partial x_j}=\beta p_i(\delta_{ij}-p_j),
\end{equation}
\begin{equation} \label{partial_max_jk}
	\;\; \frac{\partial^2 p_i}{\partial x_j \partial x_k}=\beta^2 p_i\left[(\delta_{ik}-p_k)\delta_{ij}-p_j(\delta_{ik}+\delta_{jk}-2p_k)\right].
\end{equation}

First, we fix $i\neq j$, the first term in \eqref{VarDiff_04} becomes $p_j(Z(\theta))p_i(Z(\theta))$ which is positive. In order to answer Problem \ref{problem1}, it suffices to show that the second term in \eqref{VarDiff_04}
\begin{equation} \label{Covariance_max}
	-\text{Cov}\left(Q_\beta(Z(\theta)) ,\beta p_i(Z(\theta))p_j(Z(\theta)) \right)\geq 0,
\end{equation}
for each $i\neq j, \theta$, given Condition \ref{condition1}, namely, $\sigma^Y_{i,j}\geq \sigma^X_{i,j},\; \forall i\neq j$ and $\sigma^Y_{i,i}= \sigma^X_{i,i},\; \forall i$. \eqref{Covariance_max} is the Gaussian covariance inequality we need.\\

In the following, we discuss our results with literature. In \eqref{Covariance_max}, $Q_\beta$ is increasing in each argument, while $p_ip_j$ is decreasing in $x_k$, $k\neq i,j$, separately and possibly increasing in $x_i$ or $x_j$. On another hand, $Q_\beta$ is a permutable convex function, and $p_ip_j$ is a log-concave function. Thus, \eqref{Covariance_max} is a question about proving a Gaussian covariance inequality for a specific pair of convex/log-concave functions.
Gaussian covariance inequalities hold for many symmetric (with respect to the origin) cases. \cite{Royen14} proved the long-standing conjecture for a class of probability distributions. Two symmetric convex sets are positively correlated. Or equivalently, one can replace the two symmetric convex sets by two even quasi-concave functions (symmetric quasi-concave/symmetric quasi-concave)\footnote{Two facts: a non-negative concave function is log-concave; a log-concave function is quasi-concave.}. See \cite{LatMat17} for a proof only for Gaussian measures and the reference therein.
However, in \eqref{Covariance_max}, none of $Q_\beta$ and $p_ip_j$ is symmetric about the origin.
For non-symmetric cases, the result in \cite{Harge04} needs a drift correction. 
\begin{theorem}\label{thm:Harge}  (\cite{Harge04}) Let $f$ be a convex function on $\bbR^N$ and $g$ a log-concave function on $\bbR^N$. 
	Let $\mu$ be a Gaussian measure on $\bbR^N$. Then
\begin{equation} \label{Harge-cor}
	\int f(x+l-m) \frac{g(x)d\mu(x)}{\int g d\mu} \leq \int f\; d\mu,
\end{equation}
where
\begin{equation} \label{Harge-mean}
\displaystyle l =\int x \; d\mu, \;\; m=\int x\frac{g(x)d\mu(x)}{\int g d\mu}.
\end{equation}
	
\end{theorem}

When the Gaussian measure $\mu$ is centered ($l=\vec{0}$) and the log-concave function $g$ is even, then $m=\vec{0}$. 
The drift $l-m$ goes away (non-symmetric convex/symmetric log-concave).
An non-symmetric example without a drift correction is in \cite{SW99}. The result in \cite{SW99} showed that a convex set and a strip are positively correlated under ``centroid condition'' for any Gaussian measures. ``Centroid condition'' roughly says that the two sets have the same directional bias.

We prove a specific case for a (non-symmetric convex/non-symmetric log-concave) correlation inequality for i.i.d. standard Gaussians.
Denote especialy $\bG\sim Normal(0_N, id_N)$, we have

\begin{theorem}\label{thm_iid}  $N\geq 3$, $i\neq j$, $\beta>0$,
\begin{equation}\label{inequality_iid}
	\text{Cov}(\log S_N(\bG), p_i(\bG)p_j(\bG))\leq 0.\end{equation}
\end{theorem}

\begin{remark}\label{rmk:bivariate}
	When $N=2$, \eqref{inequality_iid} is true for any bivariate Gaussian distrbution.
\end{remark}

The other main result in this paper is a partial answer to Problem \ref{problem1}, by assuming that $Y$ has only the first two standard Gaussians positively correlated.
\begin{theorem}\label{thm_rho} $N\geq 3$.
	Under Assumption \ref{assump1} and let $Y\sim Normal(0, \Sigma^Y)$ with $\sigma^Y_{i,j}=0$, for $i\neq j$, except $\sigma^Y_{1,2}=\rho \in(0,1)$. Let $X$ be the standard normal vector. Then
	$$\text{Var}(M_N(Y))\geq \text{Var}(M_N(X)) .$$
\end{theorem}

\begin{corollary}\label{thm_decreasing} $N\geq 3$,
	$$\text{Var}(M_{N-1}(\bG))\geq \text{Var}(M_{N}(\bG)) .$$
\end{corollary}

\begin{remark}
	 It is well-known that $\text{Var}(M_{N}(\bG)) \sim \frac{1}{2\log N}$, see \cite{LLR83}.
\end{remark}

\begin{remark}\label{bivariate_max}
Let $G$ be any bivariate Gaussian distribution with covariance matrix $\cC$. 
$$\text{Var} (M_2(G))=\frac{(c_{1,1}+c_{2,2})}{2}\left(1-\frac{1}{\pi}\right)+ \frac{1}{\pi}c_{1,2}.$$ Therefore, $\text{Var} (M_2(Y))=1-\frac{1}{\pi}+\frac{1}{\pi}\sigma_{1,2}^Y$ and $\text{Var} (M_2(\bG))=1-\frac{1}{\pi}$. This idicates that Theorem \ref{thm_rho} and Corollary \ref{thm_decreasing} are true for $N=2$.
\end{remark}

At the end of the introduction, we make the following conjecture.
\begin{conjecture}\label{conj1}
	Under Assumption \ref{assump1}, $\sigma^Y_{i,j}\geq \sigma^X_{i,j} \geq 0$ is a sufficient condition for \eqref{inequality_iid}.
\end{conjecture}
\noindent
If Conjecture \ref{conj1} were true, by \eqref{VarDiff_03} and \eqref{VarDiff_04}, $\sigma^Y_{i,j}\geq \sigma^X_{i,j}\geq 0,\; \forall i\neq j$ is a sufficient condition for Problem \ref{problem1}. 

The rest of paper will be organized as the following. Section 2 will provide the proofs for Remark \ref{rmk:bivariate} and Theorem \ref{thm_iid}, and Section 3 is devoted to Theorem \ref{thm_rho}. Finally, we discuss possible applications of Problem \ref{problem1} and Corollary \ref{thm_decreasing} in Section 4.

\section{Correlation inequalities for the i.i.d. case}

\medskip
We first discuss Remark \ref{rmk:bivariate}. 
Let $\cC$ be an $N\times N$ symmetric positive semi-definite matrix and $G\sim Normal(0_N,\cC)$.
For $N=2$, $S_2(G)$ can be rewritten as 
$$S_2(G) = e^{\beta G_1}+e^{\beta G_2} = e^{\beta \frac{G_1+G_2}{2}}2\cosh\left(\beta \frac{G_1-G_2}{2}\right). $$
Therefore,
$$p_1(G)p_2(G)= \frac{e^{\beta G_1+\beta G_2}}{S_N^2(G)}=\frac{1}{\left[2\cosh\left(\beta \frac{G_1-G_2}{2}\right)\right]^2}.$$
This leads to
\begin{equation} \begin{array}{rl}
		&\text{Cov}(\log S_2(G), p_1(G)p_2(G))\\
	=&	\text{Cov}\left(\beta \frac{G_1+G_2}{2},\frac{1}{\left[2\cosh\left(\beta \frac{G_1-G_2}{2}\right)\right]^2}\right)
	   +\text{Cov}\left(\log(2\cosh\left(\beta \frac{G_1-G_2}{2}\right)),\frac{1}{\left[2\cosh\left(\beta \frac{G_1-G_2}{2}\right)\right]^2}\right)\\
	\leq & \frac{\beta}{2} E\left[(G_1+G_2) \cdot\frac{1}{\left[2\cosh\left(\beta \frac{G_1-G_2}{2}\right)\right]^2}\right], 
			\end{array}
\end{equation}
since $\frac{1}{2}\log x$ is increasing and $1/x$ is decreasing.

By the conditional expection formula, $E[G_1+G_2|(G_1-G_2)]=\frac{\text{Var}(G_1)-\text{Var}(G_2)}{\text{Var}(G_1-G_2)}(G_1-G_2)$, which is an odd function of $(G_1-G_2)$. The property that $\cosh(x)$ is even gives
$$E\left[(G_1+G_2) \cdot\frac{1}{\left[2\cosh\left(\beta \frac{G_1-G_2}{2}\right)\right]^2}\right] =0.$$
We have that $\text{Cov}(\log S_2(G), p_1(G)p_2(G)) \leq 0$ for $\beta>0$ and any bivariate Gaussian distribution.

\medskip
For $N\geq3$, we need the following covariance equality. 
Let $\phi, \psi: \bbR^N \to \bbR$ be two smooth functions. Then a covariance equality gives
\begin{equation}\label{interpolation_Nourdin}
	\text{Cov}(\phi(G), \psi(G)) =\int_{0}^1 E\:\langle \sqrt{\cC}\: \nabla \phi (G_b), \sqrt{\cC} \: \nabla \psi (H_b)\rangle_{\bbR^N} \: db,
\end{equation}
where $$(G_b, H_b) \sim Normal\left(0_{2N}, 
\begin{bmatrix}
	\cC & b\cC \\
	b\cC & \cC
\end{bmatrix}\right), \; \; 0\leq b\leq 1.$$
For the proof of \eqref{interpolation_Nourdin}, please see Lemma 4.1 in \cite{Nourdin12}.
Notice that $G_0$ and $H_0$ are independent copies and $G_1=H_1$.

Here, we prove Theorem \ref{thm_iid}.

\begin{proof}
Take $i=1,j=2$, $\phi= Q_\beta$ and $\psi= -\beta p_1p_2$. $\cC$ is the identity matrix in this case.
By \eqref{interpolation_Nourdin},
\begin{equation} \begin{array}{rl}\label{Main_cor_id0}
		&\displaystyle \text{Cov}\left(Q_\beta(\bG) ,-\beta p_1(\bG)p_2(\bG) \right) \\
		= &\displaystyle \int_{0}^1 \sum_{k,l=1}^N c_{k,l} E\left[p_k (G_b) \cdot\beta^2(-p_1p_2)(\delta_{1l}+\delta_{2l}-2p_l) (H_b) \right] db\\
		= &\displaystyle \int_{0}^1  E\left[p_1 (G_b) \cdot\beta^2(-p_1p_2)(1-2p_1) (H_b) \right] db\\
		&\displaystyle + \int_{0}^1  E\left[p_2 (G_b) \cdot\beta^2(-p_1p_2)(1-2p_2) (H_b) \right] db\\
		&\displaystyle + \int_{0}^1 \sum_{k=3}^N E\left[p_k (G_b) \cdot\beta^2(-p_1p_2)(-2p_k) (H_b) \right] db\\
	\end{array}
\end{equation}
Using $p_1 =\frac{p_1+p_2}{2} + \frac{p_1-p_2}{2}$ and $p_2 =\frac{p_1+p_2}{2} - \frac{p_1-p_2}{2}$
to rearrange the second and third integrals
in \eqref{Main_cor_id0}, and in the forth integral, we change $p_k(G_b)(p_1p_2p_k)(H_b)$ to $p_1(G_b)(p_kp_2p_1)(H_b)$
because of the exchangeability of $p_1$ and $p_k$. \eqref{Main_cor_id0} becomes
\begin{equation} \begin{array}{rl}\label{Main_cor_id1}
		&\displaystyle \text{Cov}\left(Q_\beta(\bG) ,-\beta p_1(\bG)p_2(\bG) \right) \\
				= &\displaystyle \int_{0}^1  E\left[\frac{p_1(G_b)+p_2(G_b)}{2} \cdot\beta^2(-p_1p_2)(\sum_{k=3}^N 2p_k (H_b)) \right] db\\
		&\displaystyle + \int_{0}^1  E\left[\frac{p_1(G_b)-p_2(G_b)}{2} \cdot\beta^2(-p_1p_2)(2p_2-2p_1) (H_b) \right] db\\
		&\displaystyle + \int_{0}^1 \sum_{k=3}^N E\left[p_1 (G_b) \cdot\beta^2(p_1p_2)(2p_k) (H_b) \right] db \;\;(\text{because of the exchangeability})
	\end{array}
\end{equation}	
The integrand inside the first integral of \eqref{Main_cor_id1} is equal to $E\left[ p_1(G_b)\cdot\beta^2(-p_1p_2)(\sum\limits_{k=3}^N 2p_k (H_b)) \right]$, again by the exchangeability of $p_1$ and $p_2$. Thus, the first integral and the third one in \eqref{Main_cor_id1} are cancelled out. As a result, \eqref{Main_cor_id1} becomes
\begin{equation}\label{Main_cor_id2}
	\text{Cov}\left(Q_\beta(\bG) ,-\beta p_1(\bG)p_2(\bG) \right) 
= \beta^2\int_{0}^1  E\left[(p_1-p_2)(G_b) \cdot p_1p_2(p_1-p_2) (H_b) \right] db.
\end{equation}	
\noindent	
Set
\[ \begin{array}{rl}
	x^+_{1,2}&:=\frac{1}{\sqrt{2}}(x_1+x_2),  \\
	x^-_{1,2}&:=\frac{1}{\sqrt{2}}(x_1-x_2). 
\end{array}\]
Notice that fix $b\in[0,1]$, $G^-_{b,1,2}=\frac{1}{\sqrt{2}}(G_{b,1}-G_{b,2})$ is independent of $G^+_{b,1,2}=\frac{1}{\sqrt{2}}(G_{b,1}+G_{b,2})$ and $G_{b,i}, i=3,...,N$. Furthermore, we have that
$$p_1(G_b)-p_2(G_b)=\frac{e^{\beta \frac{1}{\sqrt{2}} G_{b,1,2}^+ } \cdot 2\sinh\left(\beta \frac{1}{\sqrt{2}} G_{b,1,2}^-\right)}{e^{\beta \frac{1}{\sqrt{2}} G_{b,1,2}^+ } \cdot 2\cosh\left(\beta \frac{1}{\sqrt{2}} G_{b,1,2}^-\right)+\sum_{i\neq 1,2}e^{\beta G_{b,i}}}$$
is increasing in $G^-_{b,1,2}$ and
$$p_1p_2(p_1-p_2)(H_b)= \frac{e^{\beta \frac{2}{\sqrt{2}} H_{b,1,2}^+ } \cdot e^{\beta \frac{1}{\sqrt{2}} H_{b,1,2}^+ } \cdot 2\sinh\left(\beta \frac{1}{\sqrt{2}} H_{b,1,2}^-\right)}{\left(e^{\beta \frac{1}{\sqrt{2}} H_{b,1,2}^+ } \cdot 2\cosh\left(\beta \frac{1}{\sqrt{2}} H_{b,1,2}^-\right)+\sum_{i\neq 1,2}e^{\beta H_{b,i}}\right)^3}$$
is an odd function in $H^-_{b,1,2}$.\\

Together with the fact taht $G^-_{b,1,2}$ is positively correlated with $H^-_{b,1,2}$, we prove \eqref{Main_cor_id2} is positive for all $\beta >0$.
	
\end{proof}

One may think every ``non-symmetric convex/non-symmetric log-concave'' pair gives negative covariance because many cases do. 
However, here is a counterexample.
\begin{corollary} $N\geq 3$, $\beta>0$,
	$\text{Cov}(\log S_N(\bG), p_1(\bG)(1-p_1(\bG)))\leq 0$. Thus, 
$$\text{Cov}(\log S_N(\bG), p_1^2(\bG))\geq 0.$$
\end{corollary}

$\log S_N(x)$ is a permutable convex function and $p_1^2(x)$ is log-concave. Corollary 2.1 gives different sign from Theorem \ref{thm_iid},
where $p_ip_j$ is also log-concave.

\section{Variance comparison for maxima}

\medskip
The following is the setting for this section. Let $Y\sim Normal(0_N, \Sigma^Y)$ with $\sigma^Y_{i,i}=1, i=1,...,N$, and $\sigma^Y_{i,j}=0$, of $i\neq j$, except $\sigma^Y_{1,2}=\rho \in(0,1)$. Let $X$ be the standard normal vector. For $0\leq \theta \leq 1$, let $Z(\theta)\sim Normal (0_N, \Sigma (\theta) )$, where $\Sigma (\theta)=id_N+\theta(\Sigma_Y-id_N)$. Thus, Var$(Z_i(\theta))=1$, Cov$(Z_i(\theta),Z_j(\theta))=0, \forall i\neq j$, except Cov$(Z_1(\theta),Z_2(\theta))=\theta \rho$. Set
\[ \begin{array}{rl}
	z^+_{1,2}&=\frac{1}{\sqrt{2}}(z_1+z_2),  \\
	z^-_{1,2}&=\frac{1}{\sqrt{2}}(z_1-z_2).
		\end{array}\]
So Var$(Z^+_{1,2}(\theta))=1+\theta\rho$, Var$(Z^-_{1,2}(\theta))=1-\theta\rho$ and Cov$(Z^+_{1,2}(\theta),Z^-_{1,2}(\theta))=0$. For every $0\leq \theta \leq 1$, $Z^+_{1,2}(\theta)$ is independent of $Z^-_{1,2}(\theta)$. To prove Thoerem \ref{thm_rho}, according to the discussion in Section 1,
we need to prove
$$\lim_{\beta \to \infty} \int^1_{0} d\theta \: \text{Cov}\left(Q_\beta(Z(\theta)) ,-\beta p_1(Z(\theta))p_2(Z(\theta)) \right)\geq 0.$$

Recall that 
$$e^{\beta z_1} + e^{\beta z_2} = \displaystyle e^{\beta\frac{1}{\sqrt{2}}z^+_{1,2}} \cdot 2\cosh \left(\beta\frac{1}{\sqrt{2}}z^-_{1,2}\right).$$
We do change of variables, set $w=\beta z^-_{1,2}$, then $W\sim Normal(0, \beta^{2} (1-\theta \rho))$ which is independent of $Z^+_{1,2},Z_3,...,Z_N$. We then have (we skip $\theta$ if there is no ambiguity)
\begin{equation} \label{eq:p1p2rho}
	\begin{array}{rl}
	  E\left[\beta p_1(Z(\theta))p_2(Z(\theta)) \right]
	& = E\left[ \beta\left(\frac{e^{\beta \frac{1}{\sqrt{2}}Z^+_{1,2}}}{ e^{\beta\frac{1}{\sqrt{2}}Z^+_{1,2}} \cdot 2\cosh \left(\frac{1}{\sqrt{2}} W\right) +\sum\limits_{3\leq i \leq N} e^{\beta Z_i}} \right)^2 \right] \\
	& = \int\limits^\infty_{-\infty} \: dw \: E\left[ \left(\frac{e^{\beta \frac{1}{\sqrt{2}}Z^+_{1,2}}}{ e^{\beta\frac{1}{\sqrt{2}}Z^+_{1,2}} \cdot 2\cosh \left(\frac{1}{\sqrt{2}} w\right) +\sum\limits_{3\leq i \leq N} e^{\beta Z_i}} \right)^2 \right] \\
	& \;\;\;\;\;\;\;\;\;\;\;\;\displaystyle \cdot \: \frac{1}{\sqrt{2\pi(1-\theta\rho)}  } e^{-\frac{w^2}{2\beta^2(1-\theta\rho)}}.
\end{array}\end{equation}
If we take $\beta\to \infty$, the integrand in \eqref{eq:p1p2rho} has the pointwise limit 
$$P \left[A^+_{1,2}(Z(\theta))\right] \cdot \frac{1}{(2\cosh (\frac{1}{\sqrt{2}} w))^2} \cdot \frac{1}{\sqrt{2\pi(1-\theta\rho)}  },$$ 
where $A^+_{1,2}=\{\frac{1}{\sqrt{2}}z^+_{1,2} > \max(z_3,...,z_N)\}$.
On another hand,
the integrand in \eqref{eq:p1p2rho} is bounded above by $\frac{1}{(2\cosh (\frac{1}{\sqrt{2}} w))^2} \frac{1}{\sqrt{2\pi(1-\theta\rho)}  }$ for all $w\in (-\infty,\infty)$, $\beta>0$ and $0\leq \theta\leq 1$. Simple calculations give $$\int^\infty_{-\infty} dw \frac{1}{(2\cosh (\frac{1}{\sqrt{2}} w))^2}=\frac{1}{\sqrt{2}}$$ and $$\int^1_{0} d\theta \frac{1}{\sqrt{2\pi(1-\theta\rho)}  }=\frac{1}{\sqrt{2\pi}} \cdot \frac{2}{1+\sqrt{1-\rho}}.$$ 
Therefore, the integrand in \eqref{eq:p1p2rho} has an integrable upper bound for all $\beta>0$.

Let $\beta\to \infty$, by dominate convergence theorem,
\begin{equation}\label{key_max_01}
 E\left[\beta p_1(Z(\theta))p_2(Z(\theta)) \right]\to   P \left[A^+_{1,2}\right] \cdot I
\end{equation}
and
\begin{equation}\label{key_max_01}
\int_0^1 d\theta \: E\left[\beta p_1(Z(\theta))p_2(Z(\theta)) \right]\to  \int_0^1 d\theta \: P \left[A^+_{1,2}\right] \cdot I,
\end{equation}
where
\begin{equation}
	I=I(\rho,\theta):=\frac{1}{2\sqrt{\pi (1-\theta\rho)} }.
\end{equation}
This also gives $\int_0^1 d\theta \: E\left[ p_1(Z(\theta))p_2(Z(\theta)) \right]\to 0$ as $\beta\to \infty$.

On another hand,
$$Q_\beta(z)=\frac{1}{\beta} \log \left(e^{\beta\frac{1}{\sqrt{2}}z^+_{1,2}} \cdot 2\cosh \left(\beta \frac{1}{\sqrt{2}} z^-_{1,2}\right) +\sum\limits_{3\leq i \leq N} e^{\beta z_i}\right).$$
We take $$M'_N (z):= \max(\frac{1}{\sqrt{2}}z^+_{1,2}, z_3,...,z_N)$$ 
which is less than $M_N(z)$. And

$$M'_N(z) \leq Q_\beta(z) \leq \frac{1}{\beta} \log \left( 2\cosh \left(\beta \frac{1}{\sqrt{2}} z^-_{1,2}\right) (N-1)\right) +M'_N(z)$$

\noindent
Thus, $\beta\to\infty$ gives\footnote{The integrand in this case is bounded above by another function which is integrable via elementary calculations.}
\begin{equation}\label{key_max_02}
	E\left[Q_\beta(Z(\theta)) \cdot \beta p_1(Z(\theta))p_2(Z(\theta)) \right]\to E\left[M'_N (Z(\theta)) \cdot 1_{A^+_{1,2}}\right] \cdot I.
\end{equation}
Again, $1_{A^+_{1,2}}= 1_{A^+_{1,2}}(Z(\theta))$.

As a consequence, let $\beta\to \infty$, 
\begin{equation} \label{key_max_03}
	  \text{Cov}\left(Q_\beta(Z(\theta)) ,\beta p_1(Z(\theta))p_2(Z(\theta)) \right)
	\to   \left(E\left[M'_N (Z(\theta)) \cdot 1_{A^+_{1,2}}\right] -E\left[M_N(Z(\theta))\right] P \left[A^+_{1,2}\right]\right)\cdot I
 \end{equation}
and
\begin{equation} \label{key_max_04}
	\begin{array}{rcl}
		\text{Var}(M_N(Y))-\text{Var}(M_N(X))
		  &=& \displaystyle -2 \rho \int^1_0 I(\rho,\theta) d\theta  \\
		  & & \displaystyle \cdot \left(E\left[M'_N (Z(\theta)) \cdot 1_{A^+_{1,2}}\right] -E\left[M_N(Z(\theta))\right] P \left[A^+_{1,2}\right]\right).
	\end{array}
\end{equation}
To prove Theorem \ref{thm_rho}, it is sufficient to show
\begin{equation} \label{key_max_05}
	E\left[M'_N (Z(\theta)) \cdot 1_{A^+_{1,2}}\right] \leq E\left[M'_N(Z(\theta))\right] P \left[A^+_{1,2}(Z(\theta))\right]
\end{equation}
for all $0\leq \theta \leq 1$.

In \eqref{key_max_05}, $\frac{1}{\sqrt{2}}Z^+_{1,2}(\theta)$, $Z_i(\theta)$, $i=3,...,N$ are independent Gaussians with
Var$(\frac{1}{\sqrt{2}}Z^+_{1,2}(\theta))=\frac{1+\theta\rho}{2}\leq 1$ and Var$(Z_i(\theta))=1$, $i=3,...,N$.\\

The following is another specific (non-symmetric convex/non-symmetric log-concave) correlation inequality.

\begin{lemma}\label{lem_max_01} $N\geq 2$, $G\sim Normal(0_N,\cC)$, $c_{1,1}\leq 1$, $c_{k,k}= 1$, $k=2,...,N$ and $c_{i,j}=0, \;\forall i\neq j$.
Then
\begin{equation}\label{cov_max_01}	
\text{Cov}\left(M_N(G) , 1_{A_1}(G)\right)= \lim_{\beta \to \infty}	\text{Cov}\left(Q_\beta(G) , p_1(G)\right)\leq 0,
\end{equation}	
where $A_1=\{x_1 > \max(x_2,...,x_N)\}$.
\end{lemma}
Apply Lemma \ref{lem_max_01} to \eqref{key_max_05}, Theorem \ref{thm_rho} is proved.\\

\noindent
{\it Proof of Lemma \ref{lem_max_01}.}\\

By \eqref{interpolation_Nourdin},
\begin{equation} \begin{array}{rl}
	&\text{Cov}\left(Q_\beta(G) , p_1(G)\right)\\
	&= \int_{0}^1 \: db \: \left\{   c_{1,1} E\left[  p_1 (G_b)\cdot \beta p_1(H_b)(1-p_1(H_b)) \right] -
	\sum_{k=2}^N   E\left[p_k(G_b)\cdot \beta p_1(H_b)p_k(H_b)\right]\right\}\\
	&\leq \int_{0}^1 \: db \: \left\{ \sum_{k=2}^N  E\left[(p_1(G_b)-p_k(G_b))\cdot  \beta p_1(H_b)p_k(H_b)\right] \right\}
\end{array} \end{equation}
We are going to show
\begin{equation}\label{ineq:term}
	\lim_{\beta \to \infty} E\left[(p_1(G_b)-p_k(G_b))\cdot \beta p_1(H_b)p_k(H_b)\right] \leq 0 
\end{equation}
for $k=2,...,N$.\\

Rewrite
$$p_1(G_b)-p_k(G_b)=\frac{e^{\beta \frac{1}{\sqrt{2}} G_{b,1,k}^+ } \cdot 2\sinh\left(\beta \frac{1}{\sqrt{2}} G_{b,1,k}^-\right)}{e^{\beta \frac{1}{\sqrt{2}} G_{b,1,k}^+ } \cdot 2\cosh\left(\beta \frac{1}{\sqrt{2}} G_{b,1,k}^-\right)+\sum_{i\neq 1,k}e^{\beta G_{b,i}}}$$
and
$$p_1(H_b)p_k(H_b)= \frac{e^{\beta \frac{2}{\sqrt{2}} H_{b,1,k}^+ } }{\left(e^{\beta \frac{1}{\sqrt{2}} H_{b,1,k}^+ } \cdot 2\cosh\left(\beta \frac{1}{\sqrt{2}} H_{b,1,k}^-\right)+\sum_{i\neq 1,k}e^{\beta H_{b,i}}\right)^2}.$$

$\text{Var}(G_{b,1,k}^+)=\text{Var}(G_{b,1,k}^-)=\frac{c_{1,1}+1}{2}=:\gamma_{11}$.
Because of the hyperbolic functions, we look at $$\bar{G}:=-G_{b,1,k}^-$$ and $$\bar{H}:=-H_{b,1,k}^-$$ instead of $G_{b,1,k}^-$ and $H_{b,1,k}^-$. Firstly, we have that $\langle\; G_{b,1,k}^+,\: \bar{G}, \; H_{b,1,k}^+, \; \bar{H} \;\rangle$ are independent of other $G_{b,i}$'s and $H_{b,i}$'s. We take $$\cG:=\sum_{i\neq 1,k}e^{\beta G_{b,i}}$$ and $$\cH:= \sum_{i\neq 1,k}e^{\beta H_{b,i}}.$$ 
Secondly, we denote the covariance matrix of $\langle\; G_{b,1,k}^+,\: \bar{G} \;\rangle$ as $\Gamma$ and we have $$\gamma_{12}=\text{Cov}(G_{b,1,k}^+,\bar{G})=\frac{-c_{1,1}+1}{2}\geq 0.$$ Thus,
$\langle\; G_{b,1,k}^+,\: \bar{G}, \; H_{b,1,k}^+, \; \bar{H} \;\rangle$ are positively correlated since their covariance matrix is
\[\begin{bmatrix}
	\Gamma & b\Gamma \\
	b\Gamma & \Gamma
\end{bmatrix}.\]

Represent $\langle\; G_{b,1,k}^+,\: \bar{G}, \; H_{b,1,k}^+ \;\rangle$ as their projections on the space spanned by $\bar{H}$ and the residues:
\[\begin{array}{lll}
G_{b,1,k}^+ & = \displaystyle \frac{b\gamma_{12}}{\gamma_{11}}\bar{H} + G^+_\delta & := \frac{b\gamma_{12}}{\gamma_{11}}\bar{H}  + \left\{G_{b,1,k}^+- \frac{b\gamma_{12}}{\gamma_{11}}\bar{H} \right\},\\
   \bar{G}  & = \displaystyle b \bar{H} + G^-_\delta & := b \bar{H} +\left\{\bar{G} - b \bar{H}\right\}, \\
H_{b,1,k}^+ &=  \displaystyle \frac{\gamma_{12}}{\gamma_{11}}\bar{H} + H^+_\delta & := \frac{\gamma_{12}}{\gamma_{11}}\bar{H}  + \left\{H_{b,1,k}^+ - \frac{\gamma_{12}}{\gamma_{11}}\bar{H}\right\}.
\end{array}\]

The term in \eqref{ineq:term} with a minus sign becomes
\[\begin{array}{rl}
	&-E\left[ (p_1(G_b)-p_k(G_b))\cdot p_1(H_b)p_k(H_b)\right]\\
=& E\left[ \frac{e^{\beta \frac{1}{\sqrt{2}} \left( \frac{b\gamma_{12}}{\gamma_{11}}\bar{H}  + G^+_\delta \right) } \cdot 2\sinh\left(\beta \frac{1}{\sqrt{2}} \left(b \bar{H} +G^-_\delta\right)\right)}{e^{\beta \frac{1}{\sqrt{2}} \left( \frac{b\gamma_{12}}{\gamma_{11}}\bar{H}  + G^+_\delta \right) } \cdot 2\cosh\left(\beta \frac{1}{\sqrt{2}} \left(b \bar{H} +G^-_\delta\right)\right)+ \cG}
     \cdot  \frac{e^{\beta \frac{2}{\sqrt{2}} \left(\frac{\gamma_{12}}{\gamma_{11}}\bar{H}  + H^+_\delta \right)} }{\left(e^{\beta \frac{1}{\sqrt{2}} \left(\frac{\gamma_{12}}{\gamma_{11}}\bar{H}  + H^+_\delta \right) } \cdot 2\cosh\left(\beta \frac{1}{\sqrt{2}} \bar{H}\right)+ \cH \right)^2} \right].
\end{array}
\]
It has the same limit as, if one takes $\beta\to\infty$,
\begin{equation}\label{final}
 E\left[ \frac{e^{\beta \frac{1}{\sqrt{2}} \left(   G^+_\delta \right) } \cdot 2\sinh\left(\beta \frac{1}{\sqrt{2}} \left(G^-_\delta\right)\right)}{e^{\beta \frac{1}{\sqrt{2}} \left(    G^+_\delta \right) } \cdot 2\cosh\left(\beta \frac{1}{\sqrt{2}} \left( G^-_\delta\right)\right)+ \cG}
\cdot  \left(\frac{e^{\beta \frac{1}{\sqrt{2}} \left(  H^+_\delta \right)} }{e^{\beta \frac{1}{\sqrt{2}} \left(  H^+_\delta \right) } + \cH }\right)^2 \right]\cdot I_2,
\end{equation}
where $I_2= \frac{1}{2\sqrt{\pi \text{Var}(\bar{H})}}$ and $\text{Var}(\bar{H})=\frac{c_{1,1}+1}{2}$.
We would like to show \eqref{final} is nonnegative.\\

We turn our attention to $\langle\: G^+_\delta,\: G^-_\delta,\: H^+_\delta \:\rangle$. We have
$$\text{Cov}(G^+_\delta,G^-_\delta)=(1-b^2)\gamma_{12}\geq 0,$$
$$\text{Cov}(G^+_\delta,H^+_\delta)=b\frac{\gamma_{11}^2-\gamma_{12}^2}{\gamma_{11}}\geq 0$$
 and 
$$\text{Cov}(G^-_\delta,H^+_\delta)=\text{Cov}(\bar{G} - b \bar{H} ,H^+_\delta)=\text{Cov}(\bar{G} - b \bar{H} ,H_{b,1,k}^+)=0.$$ 
The last property shows that $G^-_\delta$ and $H^+_\delta$ are indepednent.
Moreover, represent $G^+_\delta$ as the projection on the space spanned by $\langle\: G^-_\delta,\: H^+_\delta \:\rangle$ and the residue: 
$$G^+_\delta=\frac{\gamma_{12}}{\gamma_{11}}G^-_\delta +b H^+_\delta  + G^+_{\delta\delta}:=\frac{\gamma_{12}}{\gamma_{11}}G^-_\delta +b H^+_\delta+\left\{ G^+_\delta - \frac{\gamma_{12}}{\gamma_{11}}G^-_\delta -b H^+_\delta\right\}.$$

We finally split the expectation in \eqref{final} into $\{G^-_\delta >0\}$ and $\{G^-_\delta <0\}$ and use the fact that
$$\frac{e^{\beta \frac{1}{\sqrt{2}} \left(   \frac{\gamma_{12}}{\gamma_{11}}G^-_\delta +b H^+_\delta  + G^+_{\delta\delta} \right) } \cdot 2\sinh\left(\beta \frac{1}{\sqrt{2}} \left(G^-_\delta\right)\right)}{e^{\beta \frac{1}{\sqrt{2}} \left(    \frac{\gamma_{12}}{\gamma_{11}}G^-_\delta +b H^+_\delta  + G^+_{\delta\delta} \right) } \cdot 2\cosh\left(\beta \frac{1}{\sqrt{2}} \left( G^-_\delta\right)\right)+ \cG} 
> 
\frac{e^{\beta \frac{1}{\sqrt{2}} \left(   -\frac{\gamma_{12}}{\gamma_{11}}G^-_\delta +b H^+_\delta  + G^+_{\delta\delta} \right) } \cdot 2\sinh\left(\beta \frac{1}{\sqrt{2}} \left(G^-_\delta\right)\right)}{e^{\beta \frac{1}{\sqrt{2}} \left(    -\frac{\gamma_{12}}{\gamma_{11}}G^-_\delta +b H^+_\delta  + G^+_{\delta\delta} \right) } \cdot 2\cosh\left(\beta \frac{1}{\sqrt{2}} \left( G^-_\delta\right)\right)+ \cG} $$
when $G^-_\delta$ is positive.
Thus, \eqref{final} is nonnegative and the proof of Lemma \ref{lem_max_01} is complete.


\qed

\section{Discussions}

One possible application due to Problem \ref{problem1} is the {\it directed polymers in random environments} (DPRE), see more details in monograph \cite{Com17}. Denote $S=\{S_0,S_1,S_2,...,S_n,...\}$ as a $d$-dimensional simple random walk.
Let $H_n^{DP}(S)$ be the Hamiltonian with system size $n$ in ``Gaussian" random environments.
More precisely,
\begin{equation}
	H_n^{DP}(S) :=  \sum_{i=1}^n \sum_{x\in\bbZ^d} \omega(i,x)\cdot 1_{S_i=x},
\end{equation}
where the random environment $\omega$ is deﬁned on time and space $\bbN\times \bbZ^d$ and $\omega(i,x)$'s are i.i.d. standard normal random variables.
For each $n$-step path, 
\begin{equation}
	H_n^{DP}(S)\sim Normal(0,n).
\end{equation}
For two diﬀerent paths $S$ and $\hat{S}$,
\begin{equation}\label{DP_structure}
	Cov(H_n^{DP}(S),H_n^{DP}(\hat{S})) = \sum_{i=1}^n 1_{S_i=\hat{S}_i}\geq 0.
\end{equation}
The random variable $\max_{S} H_n^{DP}(S)$ is the maximum among $(2d)^n$ possible paths of the walk $S$ at time $n$.
In other words, $\max_{S} H_n^{DP}(S)$ is the maximum among $(2d)^n$ positively correlated Gaussian random variables. 
People believe that the extreme statisstic $\max_{S} H_n^{DP}(S)$ reveals some property of the system. 
Therefore, researchers concern the scale of the following variance
\begin{equation}\label{DP_variance}
Var\left(\max_{S}H_n^{DP}(S)\right).
\end{equation}
By manipulating the covariance structure \eqref{DP_structure}, Problem \ref{problem1} may give sharp upper bounds or lower bounds for \eqref{DP_variance}.

One application due to Corollary \ref{thm_decreasing} is the confidence interval. 
Let $X_1,...,X_n$ be i.i.d. normal random variables with mean zero and variance $\sigma^2$. From Theorem 1.5.3 in \cite{LLR83},
denote 
$$a_n = (2 \log n)^{-1/2}$$
$$b_n = \sqrt{2\log n} - \frac{\log\log n +\log(4\pi)}{2\sqrt{2\log n}}$$
$$M_n^{(\sigma)} = \max_{i=1,..,n} X_i$$
and $M_n(\bG) = M_n^{(1)}$ for the case $\sigma=1$. Also let $V_n =Var(M_n(\bG))\sim a_n^2$.
The limiting distribution of $M_n^{(\sigma)}$ is
the Gumble distribution with the distribution function $Gum(x) = \text{exp}(-e^{-x})$ defined on $(-\infty, \infty)$,
that is, $\frac{M_n^{(\sigma)} -\sigma b_n}{\sigma a_n}$ converges to $Gum$ in distribution.
Since $a_n b_n \sim 1$, if we take 
\begin{equation}
	\hat{\sigma}=\frac{M_n^{(\sigma)}}{\sqrt{2\log n}},
\end{equation}
$\hat{\sigma}$ is a consistent estimator with $Var(\hat{\sigma})\sim \sigma^2a_n^4$
and $E[\hat{\sigma}] \sim \sigma$. 
The confidence interval with confidence level $\alpha$ could be derived from
\begin{equation}\label{CI_statistic}
	Gum^{-1}\left(\frac{\alpha}{2}\right) \leq \frac{\hat{\sigma} -\sigma}{ \sigma a_n \sqrt{V_n}} \leq Gum^{-1}\left(1-\frac{\alpha}{2}\right).
\end{equation}
Notice that if $\alpha \leq 20\%$, $Gum^{-1}\left(\frac{\alpha}{2}\right)<0$ due to $Gum(0)=\frac{1}{e}\approx 0.37$. 
Since $a_n \sqrt{V_n}\to 0$, when $n$ is large,
\eqref{CI_statistic} is equivalent to
\begin{equation}
	\frac{\hat{\sigma}}{1+ Gum^{-1}\left(1-\frac{\alpha}{2}\right) a_n \sqrt{V_n}} \leq \sigma
	\leq \frac{\hat{\sigma}}{1+Gum^{-1}\left(\frac{\alpha}{2}\right) a_n \sqrt{V_n}}.
\end{equation}
Corollary \ref{thm_decreasing} gives that $a_n \sqrt{V_n}$ is decreasing in $n$. Therefore, the confidence interval for $\sigma$ is narrower when $n$ is larger.

\section*{Acknowledgments}

C.-H. Huang was supported by the Ministry of Science and Technology, Taiwan, grant MOST 110-2115-M-004-001.


\end{document}